\documentclass[journal]{IEEEtran}
\ifCLASSINFOpdf
\usepackage[pdftex]{graphicx}
\graphicspath{{./}{./figures/}}
\else
\fi
\newcommand*\diff{\mathop{}\!\mathrm{d}}
\usepackage[utf8]{inputenc}
\usepackage[T1]{fontenc}
\usepackage{amsmath,color}
\usepackage{amsfonts}
\usepackage{multirow}
\usepackage[utf8]{inputenc}
\usepackage[T1]{fontenc}
\usepackage{float}
\usepackage{array}

\DeclareMathOperator{\sign}{sign}

\usepackage{hyperref}
\hypersetup{
    unicode=false,          % non-Latin characters in Acrobat’s bookmarks
    pdftoolbar=true,        % show Acrobat’s toolbar?
    pdfmenubar=true,        % show Acrobat’s menu?
    pdffitwindow=false,     % window fit to page when opened
    pdfstartview={FitH},    % fits the width of the page to the window
    pdfnewwindow=true,      % links in new PDF window
    colorlinks=true,       % false: boxed links; true: colored links
    linkcolor=red,          % color of internal links (change box color with linkbordercolor)
    citecolor=blue,        % color of links to bibliography
    filecolor=magenta,      % color of file links
    urlcolor=cyan           % color of external links
}

\begin{document}

\title{Motion control for autonomous heterogeneous multi-agent area search in uncertain conditions}

\author{Stefan Ivić\thanks{Stefan Ivić is with Faculty of Engineering, University of Rijeka, Rijeka, Croatia, e-mail: {stefan.ivic@riteh.hr}.}}

% The paper headers
\markboth{Submitted to IEEE}%
{Ivić: Motion control for autonomous heterogeneous multi-agent area search in uncertain conditions}

\maketitle
\begin{abstract}
Using multiple mobile robots in search missions offers a lot of benefits, but one needs a suitable and competent motion control algorithm which is able to consider sensors characteristics, the uncertainty of target detection and complexity of needed  maneuvers in order to make a multi-agent search autonomous. This paper provides a methodology for an autonomous two-dimensional search using multiple unmanned search agents. The proposed methodology relies on an accurate calculation of target occurrence probability distribution based on the initial estimated target distribution and continuous action of spatial variant search agent sensors. The core of the autonomous search process is a high-level motion control for multiple search agents which utilizes the probabilistic model of target occurrence via Heat Equation Driven Area Coverage (HEDAC) method. This centralized motion control algorithm is tailored for handling a group of search agents which are heterogeneous in both motion and sensing characteristics. The motion of agents is directed by the gradient of the potential field which provides near-ergodic exploration of the search space. The proposed method is tested on three realistic search mission simulations and compared with three alternative methods, where HEDAC outperforms all alternatives in all tests. Conventional search strategies need about double the time to achieve proportionate detection rate when compared to HEDAC controlled search. The scalability test showed that increasing the number of HEDAC controlled search agents, although somewhat deteriorating the search efficiency, provides needed speed-up of the search. This study shows the flexibility and competence of the proposed method and gives a strong foundation for possible real-world applications.
\end{abstract}

\begin{IEEEkeywords}
area search,
uncertainty search,
multi-agent control,
target detection,
HEDAC
\end{IEEEkeywords}

%\end{frontmatter}

\section{Introduction}

The task of performing an autonomous multi-agent search in uncertain conditions is a thriving concept which has not yet been scientifically entirely answered. Consequently, the real-world applications are not as straightforward as one would expect given the rapid development of the associated equipment.
Possible applications of two-dimensional search are numerous: search and rescue missions on land or sea, objects search and detection, surveillance, monitoring in agriculture etc. Performing a search with unmanned remote-controlled vehicles, whether terrestrial, aerial, sea surface or even submerged, gives an opportunity to organize smart and efficient motion of search agents. The use of multiple autonomous agents in the search raises a problem of coordination, which can potentially be more troublesome due to differences in motion or sensing of used search agents.
The design of reliable, robust and efficient motion control of multiple heterogeneous agents for a search in uncertain conditions is a challenging task, solution of which would provide foundations for successful real-world applications.

A general fundamental and theoretical study on the target search in uncertainty is presented in Koopman's "The Theory of Search" triad \cite{koopman1956theory1,koopman1956theory2,koopman1957theory3}. The first paper covers the kinematic bases of the search theory, specifically the motion of search agents and their sensing range, and the probability distribution of targets. 
The second sequel covers a target detection under uncertainty, with established mathematical fundamentals of spatial and temporal detection probability for both stationary and moving search agents and targets.
In the third part the search is considered as the optimization problem formulated as finding agent trajectories
which maximize the overall probability of detection for a given target occurrence probability distribution. 

There are many different approaches on high-level multi-agent search control, but two rather prominent directions stand out in the literature covering this field. These two are the approaches based on Receding Horizon Control and Spectral Multiscale Coverage methods.

\subsection{Spectral Multiscale Coverage approach}

Spectral Multiscale Coverage (SMC) was first introduced in \cite{mathew2009spectral} as an algorithm for ensuring given coverage density. 
SMC employs a gradient-based agent regulation which relies on smoothed difference between achieved and goal coverage fields. The smoothing is achieved using spectral techniques, i.e. suitably adapting the Fourier transform coefficients in order to emphasize lower modes.

The SMC idea is adopted for the 2D search in the method called Multiscale Adaptive Search presented in \cite{hubenko2011multiscale} and uses SMC on logarithmic prior in order to encompass the exponential detection model. 
Another approach, in \cite{surana2012coverage}, also employs SMC method for continuous-space uncertainty area search. Here the Sequential Ratio Probability Test and Recursive Least Squares estimation are used in order to quantify the current uncertainty in target detection and location, respectively.
The papers \cite{mathew2013experimental} and \cite{shilov2018efficient} show the experimental validation of SMC/MAS algorithm using unmanned aerial vehicles (UAV's).
The consideration of unsteady coverage problem, presented in \cite{mathew2010uniform}, is a very interesting approach and allows many useful enhancements of the method or its applications, such as the search for moving targets.
SMC method allows relatively easy implementation of kinematic and dynamic \cite{mathew2009spectral} as well as Dubins motion models \cite{surana2012coverage,mathew2013experimental}, and it is suitable to adapt to heterogeneous multi-agent control (although such adaptation have not been made yet, according to the available literature). The main drawbacks of the SMC method are the limitation to rectangular domains and lack of regulation mechanism between global and local search behavior.

\subsection{Receding Horizon Control approach}	
	
In addition to SMC based approaches, another commonly used method in coverage problems is Receding Horizon Control (RHC), also known as Model Predictive Control (MPC). The RHC method transforms the problem of discovering ergodic trajectories to the optimization problem, where agent heading directions are determined by optimizing forthcoming trajectories according to a given objective. This brings both benefits and detriments to the search procedure. The optimization procedure can easily include additional objectives or constraints such as motion model constraints, obstacle avoidance, minimizing fuel consumption etc. However, the effect of the chosen horizon, i.e. the length of trajectories which is subjected to the optimization, is twofold. A short horizon localizes the search and weakness the ergodicity of achieved trajectories but facilitates the optimization. On the other hand, long horizon enables global ergodic exploration while hindering the optimization and the efficiency.

A general coverage problem, where uniform goal coverage is considered in the domain with an obstacle, is presented in \cite{ahmadzadeh2006cooperative,ahmadzadeh2007cooperative}.
Ergodic multi-agent coverage using RHC method, with a goal of achieving a more general, non-uniform, goal coverage density is considered in \cite{mavrommati2017real,ayvali2017ergodic}.
Receding Horizon Control is utilized for uncertainty search in \cite{gan2011multi} and \cite{lanillos2014multi-uav} ensuring reduction of time needed for target detection. Arbitrary target probability maps and sensor detection models, as well as multi-agent coordination are built-in features of the control method.
A control method, designed for exploration of unknown environments while minimizing tracked landmark and agent localization uncertainty, is used for mapping planning strategy and uncertainty-aware exploration in \cite{papachristos2017uncertainty}.
As presented in \cite{yao2017gaussian}, optimal trajectories of multiple UAV's for target search with known target distribution are assured using 
Receding Horizon Control and a Gaussian Mixture Model.
Various useful constraints such as obstacle and collision avoidance, and the simultaneous arrival of multiple UAV's to a given destination are included in the proposed method. 
A target search using the distributed RHC method combines a UAV motion planning layer and a network topology control layer, using the distributed RHC optimization, is proposed in \cite{di2015potential}.
A more sophisticated tasks, such as surface inspection in three dimensions with a multi-rotor platform \cite{bircher2016receding}, can also be accomplished with Receding Horizon Control.

\subsection{Other approaches and methods}

Engaging a Bayesian approach into target search problems has been investigated by several authors \cite{bourgault2006optimal,lin2010bayesian,xiao2019cooperative}.
In these papers, various features of the search are considered such as static or mobile target search, different search vehicle kinematics and sensor detection models.
Proposed algorithms maximize the probability of detecting a target in a fixed time window.

A cooperative framework for managing a multi-agent search system is proposed in \cite{polycarpou2001cooperativesearch,polycarpou2001cooperativecontrol}. The algorithm relies on communication between agents and environmental perception using on-board sensors in order to achieve cooperation and successful obstacles avoidance. Many real motion limitations, such as maneuverability limitations or fuel/time constraints, are included in the algorithm.
A similar control, which also considers the motion and sensing limitation, is presented in \cite{hu2013multiagent}. In order to provide coverage and topology control, the method uses Bayesian update of probability map, individually for each agent.
The optimal search in terms of maximal detection rate is realized using the Mode Good Ratio heuristic in \cite{lin2014hierarchical}. Although the algorithm takes advantage of the target occurrence probability estimation and produces a successful search, it does not coordinate multiple agents. 

On a more practical side, the use of UAV's as search agents is well recognized and probably most suitable technology for carrying out the search task. Searching for targets from the air allows a fast exploration of a terrain, due to absence of physical obstacles and responsive dynamics of UAV's.
%A fast search is very important since early detection is crucial for success of search missions. Although it also depends on other conditions, in search and rescue the survival probability notably drops as time goes by \cite{adams2007search}.
AN optimal multi-agent target search using Ant Colony Optimization (ACO) method is proposed in \cite{perez2018ant}. Similar to RHC, it  balances between the computational requirements and the efficiency of found search strategy.
Simple applications of UAV's in a search employ a single unmanned vehicle \cite{goodrich2007using,goodrich2008supporting,goodrich2009towards} or performing the search using very simple agent paths \cite{vincent2004framework}.
Another approach using single UAV, presented in \cite{carpin2013variable}, utilizes probabilistic quadtree structures for directing aerial search. 
Ground robots dynamic control and coordination for wilderness search and rescue in \cite{macwan2015multirobot} achieves initial and time-optimal search agent trajectories with ability to perform path re-planning whenever it is necessary.

\subsection{Overview of the manuscript and contribution}

Although the addressed approaches are to some extent prosperous in the area search problems, there are many drawbacks and hindrances for complete and successful implementation of autonomous search. Many of the proposed methods are complicated, inflexible, inefficient or inaccurate in at least one essential aspect of area search such as motion control, search agent dynamics model or target and detection probabilistic model. The search process proposed in this paper aims to overcome these issues by adapting a relatively new are coverage control algorithm and pairing it with accurate detection and target distribution estimation model.

A multi-agent coverage control called Heat Equation Driven Area Coverage (HEDAC) algorithm is proposed in \cite{ivic2016ergodicity}. This ergodicity-based feedback method drives agents using gradient of appropriately designed potential field in order to achieve given goal coverage density.
The agent trajectories are near-ergodic as they, in infinite time, minimize the ergodicity measure proposed in ~\cite{mathew2011metrics}.
Tests showed suitability of HEDAC method to perform control of cooperative multi-agent motion for arbitrary given goal coverage density fields.
Furthermore, HEDAC has been successfully adopted for governing multi-agent nonuniform spraying in \cite{ivic2019autonomous}.

In the proposed methodology for heterogeneous multi-agent search in uncertain conditions, the HEDAC method is tailored and accompanied with specific target search apparatus. In this paragraph we briefly present essential elements required for assembling the integrated search process, which are key contributions of this paper.
The target occurrence probability is exactly calculated considering the initial target occurrence estimation and conducted search agents coverage.
Search agent sensors characteristics, such as spatial and temporal scope, and motion properties are used as a basis in exact and mathematically rigorous establishment of actual target occurrence probability.
Search agents are directed by the calculated target occurrence distribution appointed as a source of potential in the HEDAC control method. This allows ergodic exploration of target distribution and provides near-optimal search. 
The method allows the use of heterogeneous agents, in motion and sensing characteristics, for performing the search.
The proposed HEDAC heterogeneous multi-agent search control is evaluated on three realistic search scenarios and compared with alternative algorithms - one conventional and two state-of-the-art methods. These methods are suitably modified for handling multiple heterogeneous search agents with arbitrary sensing functions.
The computational efficiency of all four methods is benchmarked on different search scenarios.
Finally, the scalability of the proposed multi-agent search is investigated showing that HEDAC can successfully control the swarm of search agents without crucially deteriorating the search efficiency.

The manuscript contains literature overview, theoretical setup of the methodology, numerical validation on three test cases, computational efficiency benchmark and scalability analysis. Animated search simulations, and details about considered search scenarios and numerical implementation of control algorithms are given in supplementary material.

\section{Search in uncertant environment}

Search agents with sensing equipment dynamically explore the search space with a goal to detect sought target, with a known estimated probability distribution of target occurrence. The detection of a target depends on both temporal and spatial characteristics of sensors. 
A 2D search problem is considered, since most of the target search problems are, or can be reduced to, two-dimensional. In this section, a general probabilistic model of target detection is presented, based on arbitrary search agents paths and sensor models.

\subsection{Target detection probability}
\label{subsec:target_detection_probability}

According to \cite{koopman1956theory2}, in the case of continuous search, $\gamma\diff t$ is the probability of target detection in a short time interval of length $\diff t$, where $\gamma \geq 0$ is called the instantaneous probability density (of detection). When the search is performed under unchanging conditions, continuously until time $t$, the probability  of detection $p(t)$ is given by
\begin{equation}
 p(t)=1-e^{-\gamma t}, \; t \geq 0.
 \label{eq:detection_probability_koopman}
\end{equation}
This formulation is valid only for a single point (target) that is continuously observed, where the detection probability approaches to 1 as time goes by.  

If $\gamma$ is considered as variable in time, $\gamma(t) \geq 0$, the probability of detection until time $t$ is
\begin{equation}
 p(t) = 1 - \exp\left({-\int_0^t \gamma(\tau) \diff \tau}\right).
\end{equation}

% {\color{cyan}
% FROM SPRAYING PAPER:

% Each agent has its spraying pattern, depending on the spray nozzle characteristics and arrangement (Figure~\ref{fig:spraying_pattern}). The spraying pattern can be defined as instantaneous spraying spatial density $\phi(\mathbf{x}_{loc})$, which is defined in local agent coordinate system $\mathbf{x}_{loc} = (\chi, \zeta)$, where $\zeta$ is heading axis and $\chi$ is axis transversal to the movement of agent with origin in the location of agent $\mathbf{z}_i$. Using translation and rotation transformation, one can easily calculate local coordinates $\mathbf{x}_{loc} = (\chi, \zeta)$ from global coordinates $\mathbf{x} \in \Omega$, according to trajectory $\mathbf{z}$ and heading $\theta$ at time $t$:
% \begin{equation}
% \mathbf{x}_{loc,i}(t) = (\mathbf{z}_i(t) - \mathbf{x}) \cdot \mathbf{R}(\theta_i(t))
% \end{equation}
% %\begin{equation}
% %\phi_i \Big( (\mathbf{z}_i(t) - \mathbf{x}) \cdot \mathbf{R}(\theta_i(t)) \Big)
% %\end{equation}
% where $\mathbf{R}$ is the rotation matrix defined with agent heading angle $\theta$ as:
% \begin{equation}
% \mathbf{R}(\theta)={\begin{bmatrix}\cos \theta &-\sin \theta \\\sin \theta &\cos \theta \\\end{bmatrix}}.
% \end{equation}
% }

Although the locations of targets are unknown to the search control method, with known search agent trajectories, the detection probability for any point $\mathbf{x}$ in the domain $\Omega \subset \mathbb{R}^2$ can be defined. For now, trajectories and directions of mobile agents, and their sensing actions accordingly, are considered as known. Realized trajectories are denoted as $\mathbf{z}_i : [0,t ] \to \mathbb{R}^2$, for $i = 1,2,\ldots, N$ where $N$ is the number of search agents. Heading direction of $i$-th search agent is defined by the direction angle $\theta_i(t): [0,t] \to [0, 2\pi]$.

\begin{figure}[!ht]
	\centering
	\includegraphics[width=\linewidth]{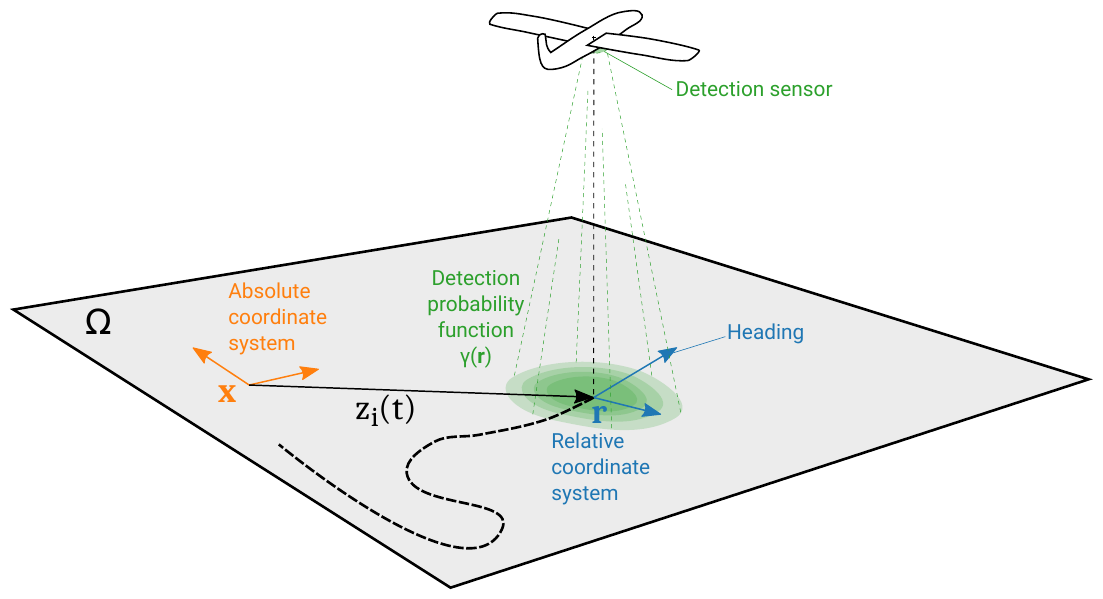}
	\caption{Absolute and relative coordinate systems, and detection function $\gamma$. (color figure available online)}
	\label{fig:sensing_function}
\end{figure}

The instantaneous detection probability can be extended to depend on the position of the target relative to the sensor, as it is the case for almost any real application. A position in agent's relative coordinate system $\mathbf{r} \in \mathbb{R}^2$ is used in order to define sensor detection probability function $\gamma = \gamma(\mathbf{r})$ (Figure~\ref{fig:sensing_function}). With origin at $\mathbf{z}_i$, relative coordinates are easily calculated as
\begin{equation}
\mathbf{r}_{i}(t) = (\mathbf{z}_i(t) - \mathbf{x}) \cdot \mathbf{R}(\theta_i(t))
\end{equation}
%\begin{equation}
%\phi_i \Big( (\mathbf{z}_i(t) - \mathbf{x}) \cdot \mathbf{R}(\theta_i(t)) \Big)
%\end{equation}
where $\mathbf{R}$ is the rotation matrix which rotates the local coordinate system according to heading direction angle $\theta$:
\begin{equation}
\mathbf{R}(\theta)={\begin{bmatrix}\cos \theta &-\sin \theta \\\sin \theta &\cos \theta \\\end{bmatrix}}.
\end{equation}

For a given sensor detection function $\gamma(\mathbf{r})$ one can easily calculate the sensing intensity of the detection sensor as
\begin{equation}
	I = \int_{\Omega_\gamma} \gamma(\mathbf{r}) \diff\mathbf{r}
\end{equation}
where $\Omega_\gamma$ is the scope of the sensor considering its maximal reach. The sensing intensity $I$ combines the spatial scope and magnitude of sensor function $\gamma$ in order to provide unambiguously measure of sensor capability.

Keeping in mind a single sensor detection function $\gamma(\mathbf{r})$, the probability of detecting a certain target at location $\mathbf{x}$ until time $t$ can be written as:
\begin{equation}
p(\mathbf{x}, t) = 1 - \exp\left({-\int_0^t \gamma(\mathbf{r}(\tau)) \diff \tau}\right)
\label{eq:point_probability}
\end{equation} 
where $\mathbf{r}(\tau)$ is relative position between agent and the target at time $\tau$.

Now, the spatial and temporal detection probability function which considers sensory activity from multiple search agents can be defined:
\begin{equation}
	p(\mathbf{x}, t) =  1 - \exp\left(-\sum_{i=1}^N \int_0^t  \gamma_i(\mathbf{r}_i (\tau )) \diff \tau \right)
	\label{eq:xt_detection_probabiliy}
\end{equation}
where the term $\gamma_i\left(\mathbf{r}_{i}(\tau)\right) = \gamma_i\left(\mathbf{z}_i(t) - \mathbf{x}) \cdot \mathbf{R}(\theta_i(t)\right)$ gives the instantaneous detection probability at location $\mathbf{x}$ affected by the $i$-th agent location $\mathbf{z}_i(\tau)$ and direction $\theta_i(\tau)$ and spatial sensor scope $\gamma_i$.

For simplicity, the coverage density field is introduced which indicates the influence of search agent sensors. The coverage field represents the accumulation of sensor coverages of all agents and it is defined as follows:
\begin{equation}
	c(\mathbf{x},t) = \sum_{i=1}^N \int_0^t  \gamma_i(\mathbf{r}_i (\tau )) \diff \tau.
	\label{eq:coverage}
\end{equation}
Using the coverage $c$, the expression \eqref{eq:xt_detection_probabiliy} can be written as
\begin{equation}
	p(\mathbf{x},t) = 1 - e^{- c(\mathbf{x}, t)}.
	\label{eq:spatial_detection_probability}
\end{equation}

\subsection{Target occurrence}

Initial (at the start of the search, $t=0$) spatial probability density of target occurrence $m_0(\mathbf{x})$ is considered as known. This is a spatial field which, for real-world applications, can be estimated or calculated based on available information, and should respect normalization:
\begin{equation}
\int_\Omega m_0(\mathbf{x}) \diff \mathbf{x} = 1.
\label{eq:occurence_normalization}
\end{equation}
If available $m_0$ does not conform to above criteria, it can be simply scaled  in order to achieve this.

The estimation of undetected target presence probability depends on both initial target expectation and the search agent exploration. This is a complement event of target not being present at all and/or being already detected. 
One needs to combine the target's initial occurrence estimation $m_0(\mathbf{x})$ and detection probability density $p(\mathbf{x})$, in order to calculate the undetected targets probability 
\begin{equation}
	m(\mathbf{x},t) = m_0(\mathbf{x}) \cdot (1  - p(\mathbf{x},t))
\end{equation}
or simplified and expressed with use of coverage density
\begin{equation}
	m(\mathbf{x},t) = m_0(\mathbf{x}) \cdot e^{-c(\mathbf{x},t)}.
	\label{eq:target_probability}
\end{equation}

It is rather easy to calculate the probability of target presence in the entire domain:
\begin{equation}
 E(t) = \int_\Omega m(\mathbf{x}, t) \diff \mathbf{x}.
 \label{eq:error}
\end{equation}
The total target presence probability $E(t)$ is a measure of target search success, and it is used as an error estimate for evaluating the search methods.

\section{Search control with HEDAC}
\label{sec:hedac}

In this section, mathematical models of agent motion, HEDAC control method and target detection mechanism are defined.
Details on the implementation of search simulations can be found in Supplementary materials.

\subsection{Motion models}

The control of agent motion is carried out by appointing aspired heading direction for two motion models: kinematic and Dubins. Both assume constant velocity of search agent which is suitable for legitimate search performance comparison between the two motion models.

For the formulation of motion models, the potential (temperature) field $u(\mathbf{x},t)$ is assumed as known. 
The procedure for determining potential field $u$, by using stationary heat equation, is explained in the next subsection.
For the simpler notation, let us introduce the (unit) vector field $\mathbf{u}(\mathbf{x},t)$ which is equal to normalized gradient of the potential field $u(\mathbf{x},t)$:
\begin{equation}
\mathbf{u}(\mathbf{x}, t) = \frac{\nabla u (\mathbf{x},t)}{\left|\nabla u (\mathbf{x},t)\right|}.
\end{equation} 

The kinematic motion model presupposes the movement with constant velocity but with the ability to make instantaneous changes of agent direction. Such motion is characteristic of multi-rotor drones. The model is based on simple first order differential equation:
\begin{equation}
\frac{\diff\mathbf{z}_i(t)}{ \diff t}=v_{a,i} \cdot \mathbf{u}(\mathbf{z}_i, t) ,\quad i=1,\ldots N
\label{eq:kinematic_model}
\end{equation}
with initial conditions
\begin{equation}
\mathbf{z}_i(0)=\mathbf{z}_{i,0},\quad \mathbf{z}_{i,0}\in\Omega,\quad i=1,\ldots N 
\label{eq:initial_condition}
\end{equation}
where $v_{a}$ is (constant) agent velocity magnitude.

A Dubins model, also known as unicycle model, is also considering planar motion with constant velocity, but with limited turning radius. This motion is suitable for simulating fixed-wing drones motion and it is already successfully used with HEDAC motion control in \cite{ivic2019autonomous}. Dubins model is defined with system of first order differential equations:
\begin{equation}
\begin{aligned}
\frac{\diff\mathbf{z}_i(t)}{ \diff t} & =
\left[
\begin{array}{c}
v_{a,i} \cdot \cos\theta_i(t) \\
v_{a,i} \cdot \sin\theta_i(t)
\end{array}
\right], & i=1,\ldots N\\
\frac{\diff \theta_i(t)}{\diff t} & =  \sign(\omega_i) \cdot \min (\left|\omega_i\right|, \omega_{max,i}), & i=1,\ldots N
\end{aligned}
\label{eq:dubins_equation}
\end{equation}
where $\omega_i$ and $\omega_{max,i}$ are aspired turning angular velocity and maximal turning angular velocity, respectively. The aspired change of the direction (aspired by the control method, but not necessary achievable due to turning constraint) is governed as angle of vector $\mathbf{u}(\mathbf{z}_i)$ relative to current direction $\mathbf{v}_i$ and can be easily calculated as follows:
\begin{equation}
\omega_i = \arctan\left(\frac{\mathbf{v}_i \bullet \mathbf{u}(\mathbf{z}_i,t)}{||\mathbf{v}_i||\cdot ||\mathbf{u}(\mathbf{z}_i,t)||}\right),\quad i=1,\ldots N
\label{eq:dubins_angular_velocity}
\end{equation}
where $\mathbf{v}_i$ is $i$-th agent current direction vector, while $\mathbf{v}_i \bullet \mathbf{u}$ denotes dot product of vectors $\mathbf{v}_i$ and $\mathbf{u}(\mathbf{z}_i,t)$. The minimal turning radius limitation is achieved with right hand side of second equation in system \eqref{eq:dubins_equation}. This constraint can be written in more intuitive manner as the limitation of direction change:
\begin{equation}
\left|\frac{\diff \theta_i}{\diff t}\right| = \left|\omega_i\right| \leq \omega_{max,i}
\label{eq:dubins_constraint},\quad i=1,\ldots N.
\end{equation}

For motion with constant speed $v_{a,i}$, there is a simple correlation between maximal angular velocity $\omega_{max,i}$ and minimal turning radius $R_{T,i}$:
$\omega_{max,i} = \frac{v_{a,i}}{R_{T,i}}$.

Beside the agents' initial position, the Dubins model requires initial orientation of agents. Thus we prescribe the initial condition as
$\theta_i(0) = \theta_{i,0}$, 
where $\theta_{i,0}$ is initial heading direction of $i$-th agent.

Both described motion models allow the use of heterogeneous agent's dynamics, and all four motion control methods considered in this paper are adapted to direct the motion of multiple heterogeneous agents.

\subsection{Obtaining potential field using heat equation}

In the original formulation of HEDAC method, a difference between normed goal and achieved coverage density is used as a source of heat equation \cite{ivic2016ergodicity}. The negative values of the coverage difference are trimmed in order to prevent repulsion of over-covered regions. The gradient of resulting temperature field leads the agents so that, indirectly, by acting with certain coverage, the minimization of source is achieved. 

Here, since we want to minimize it, the undetected target density $m(\mathbf{x}, t)$ is utilized as a heat source in the stationary heat equation for $u(\mathbf{x},t)$.
Since the field $u$ is smooth, we can employ its gradient to direct the agents in order to accomplish the minimization of potential $u$ which, in turn, leads to minimization of both the spatial and total probability of undetected targets occurrence ($m$ and $E$, respectively). 

At any time $t$, the field $u$ is obtained as a solution to the stationary heat equation defined as a partial differential equation 
\begin{equation}
\alpha \cdot \Delta u (\mathbf{x},t) = \beta \cdot u(\mathbf{x},t) -m(\mathbf{x},t)  
\label{eq:heat_equation}
\end{equation}
with the boundary condition
\begin{equation}
\frac{\partial u}{\partial \mathbf{n}} = 0,\textrm{ on }\partial \Omega.
\label{eq:neumann_bc}
\end{equation}
In the equation \eqref{eq:heat_equation}, $\Delta$ is a Laplace operator, while scalars $\alpha>0$ and $\beta>0$ are tunable parameters of HEDAC method. In the boundary condition $\mathbf{n}$ denotes the outward normal to the domain boundary. As the presented heat equation formulation showed good collaborative behavior without agent collision, the agent cooling used in \cite{ivic2016ergodicity} and the collision avoidance mechanism tested in \cite{ivic2019autonomous} are omitted for this application. 

\subsection{Target detection mechanism}

Considering there are many available sensors on the market, each with different spatial and temporal detection behavior, there is no unique model to describe target detection.
HEDAC high-level multi-agent motion control is able to handle arbitrary sensing models, by definig sensor function $\gamma(\mathbf{r})$, as described in section~\ref{subsec:target_detection_probability}.

Analogous to theoretical considerations in section \ref{subsec:target_detection_probability}, in a short finite time interval $\Delta t$, the probability of detecting $j$-th target, by at least one agent, depends on vicinity and sensing functions of all agents:
\begin{equation}
P_j(t, t+\Delta t) = 1 - \prod_i^N \left(e^{-\gamma_i(\mathbf{r}_{ij}(t)) \Delta t}\right)
\label{eq:target_detection_probability}
\end{equation} 
where $\mathbf{r}_{ij}$ is position of $j$-th target relative to the $i$-th search agent.
If $j$-th target position is denoted as $\mathbf{y}_j$, the agent-target relative position is
$\mathbf{r}_{ij} = (\mathbf{z_i} - \mathbf{y}_j) \cdot \mathbf{R}(\theta_i)$.

The sensor function $\gamma_i$ is considered as given and it allows the definition of arbitrary sensors. The use of different sensors, together with agent motion parameters, characterizes a fully heterogeneous multi-agent search system.

The detected targets are recorded during a search simulation which allows the calculation of detection rate. Detection rate can be calculated at any time of the search as a simple ratio:
$D(t) = \frac{n_D(t)}{n}$,
where $n_D$ is the number of detected targets and $n$ is the number of all targets.

% {\color{red}
% As the sensors reach is considered as constant, the instantaneous probability density (of target detection) can be written as follows:
% \begin{equation}
% \gamma_i(r)= 
% \begin{cases}
% \Gamma_i , & \text{ if } r < R_{D,i}\\
% 0 , 	   & \text{ otherwise,}
% \end{cases}
% \end{equation}
% where $\Gamma_i$ is the instantaneous detection rate within the radius of detection $R_{D,i}$ for $i$-th agent. Both $\Gamma$ and $R_D$ are given parameters of detection sensor, which allows the use of heterogeneous sensors. The use of different sensors, together with agent motion parameters, characterizes a fully heterogeneous multi-agent search system.

% The detected targets are recorded during a search simulation which allows the calculation of detection rate. Detection rate can be calculated at any time of the search as a simple ratio:
% $D(t) = \frac{n_D(t)}{n}$,
% where $n_D$ is the number of detected targets and $n$ is the number of all targets.}

\section{Simulation results and comparisons}

The evaluation of the proposed heterogeneous multi-agent search system is performed through detailed testing and comparison with alternative methods on three target search cases. Lawnmower method represents conventional approach while SMC and RHC are state-of-the-art multi-agent search control methods used in recent scientific publications. 
The implementations of SMC and RHC are closely related to algorithms presented in \cite{hubenko2011multiscale} and \cite{ahmadzadeh2006cooperative,ahmadzadeh2007cooperative}, respectively. Both methods are modified as it was necessary so they fit into presented search framework and to make them rightfully
comparable with HEDAC. RHC is using Particle Swarm Optimization (PSO) as underlaying optimization method, since it showed better results when compared to several optimization methods considered in other papers. It should be noted that RHC implementation is stochastic since PSO relies on random variables. 
All three alternative algorithms and their implementation are described in more details in supplementary material.

The same evaluation methodology is used for all tests.
The search simulation conduction and results interpretation is presented from two points of view: the analysis of trajectories obtained in a single-run search simulation for each scenario and the evaluation of the search performance based on Monte-Carlo simulations.

The search is performed on rectangular domains discretized with uniform orthogonal grid.% The spatial resolution of the numerical mesh is suited to the initial target density $m_0(\mathbf{x})$ and to the search agents parameters.
%The initial target occurrence density function $m_0(\mathbf{x})$ is considered as known and it is given for each case. 
In each search simulation run, 1000 targets are randomly positioned inside the search domain respecting the target distribution $m_0(\mathbf{x})$.
In order to demonstrate HEDAC's ability to direct multiple heterogeneous agents, different agent motion and sensing characteristics are used (Table~\ref{tab:agent_parameters}) within each test case. More details about each search scenario can be found in supplementary materials.

\begin{table*}[!ht]\small
	\caption{Summary of search agent motion and sensing parameters for all three test cases. In the sensing function plots, $r_1$ points left, $r_2$ points up and axis ticks  plots indicate 10 m distance. Precise detection probability functions are given in supplementary material. (color table available online)}
	\label{tab:agent_parameters}
	\begin{tabular}{ m{3cm}
					>{\centering\arraybackslash}m{2cm}
					>{\centering\arraybackslash}m{2cm}
					>{\centering\arraybackslash}m{2cm}
					>{\centering\arraybackslash}m{2cm}
					>{\centering\arraybackslash}m{2cm}
					>{\centering\arraybackslash}m{2cm}}
		\hline
		Test & \multicolumn{1}{|c|}{Test 1} & \multicolumn{3}{|c|}{Test 2} & \multicolumn{2}{|c}{Test 3} \\ 
		\hline
		Number of agents & 5 & 2 & 2 & 2 & 3 & 2 \\    
		Agent velocity $v_a$ [m/s] & 20 & 16 & 20 & 31 & 20 & 34 \\  
		Minimal turning radius $R_M$ [m] & 30 & 26 & 29 & 43 & 36 & 48 \\   
		Detection probability $\gamma$ \newline\newline \includegraphics[width=3cm]{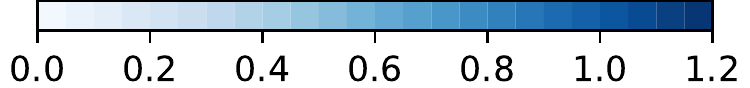}
							& \includegraphics[width=1.6cm]{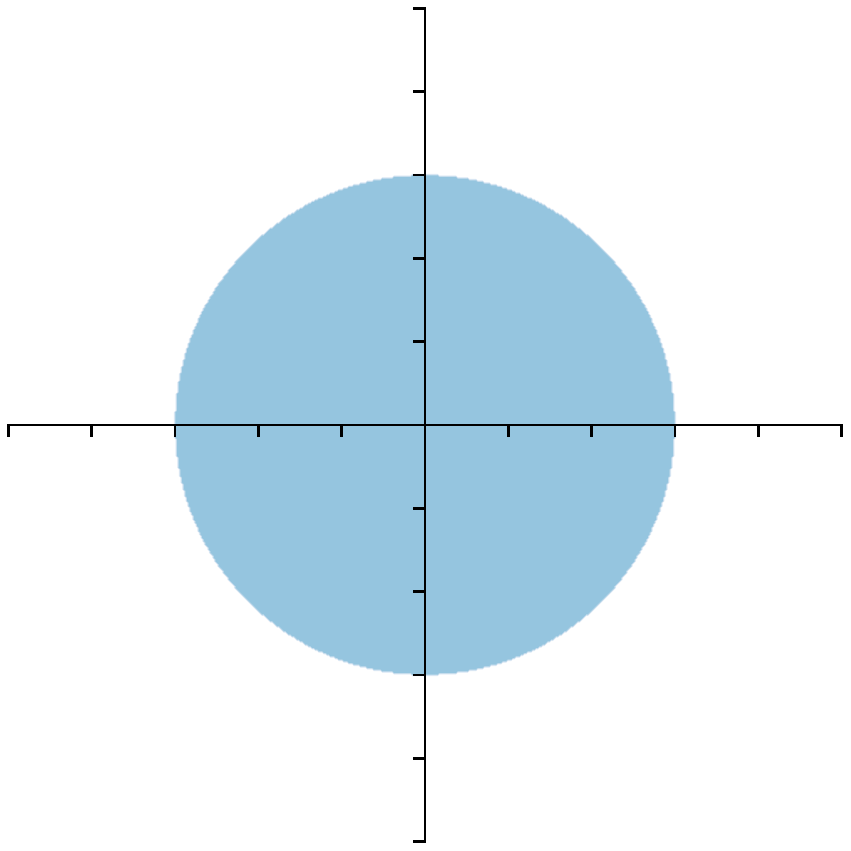}
							& \includegraphics[width=1.6cm]{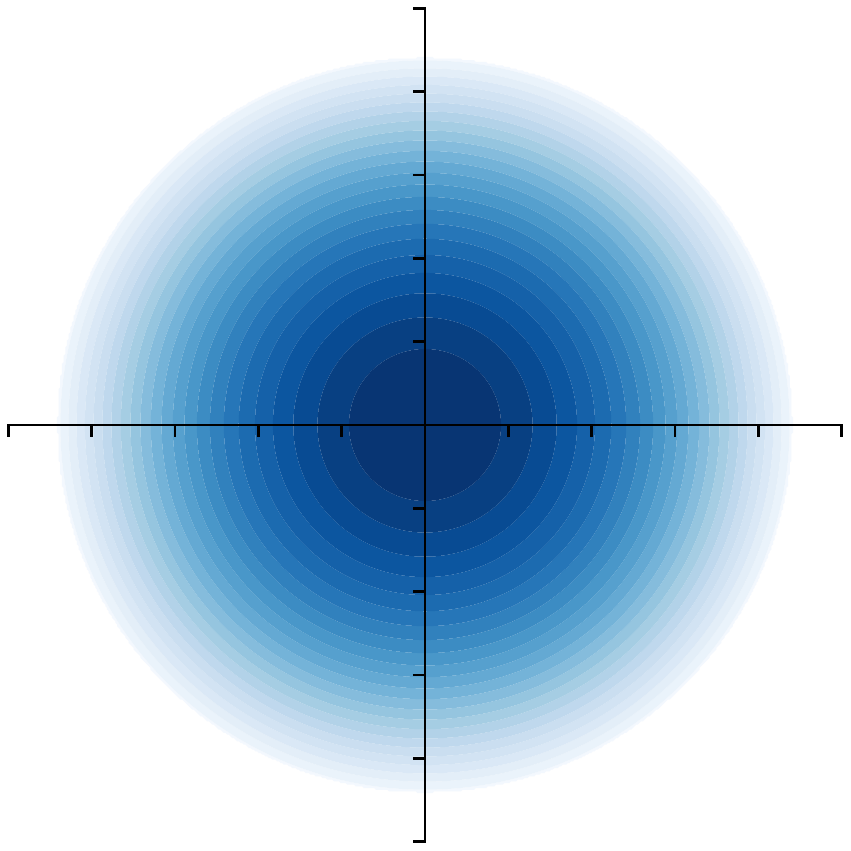}
							& \includegraphics[width=1.6cm]{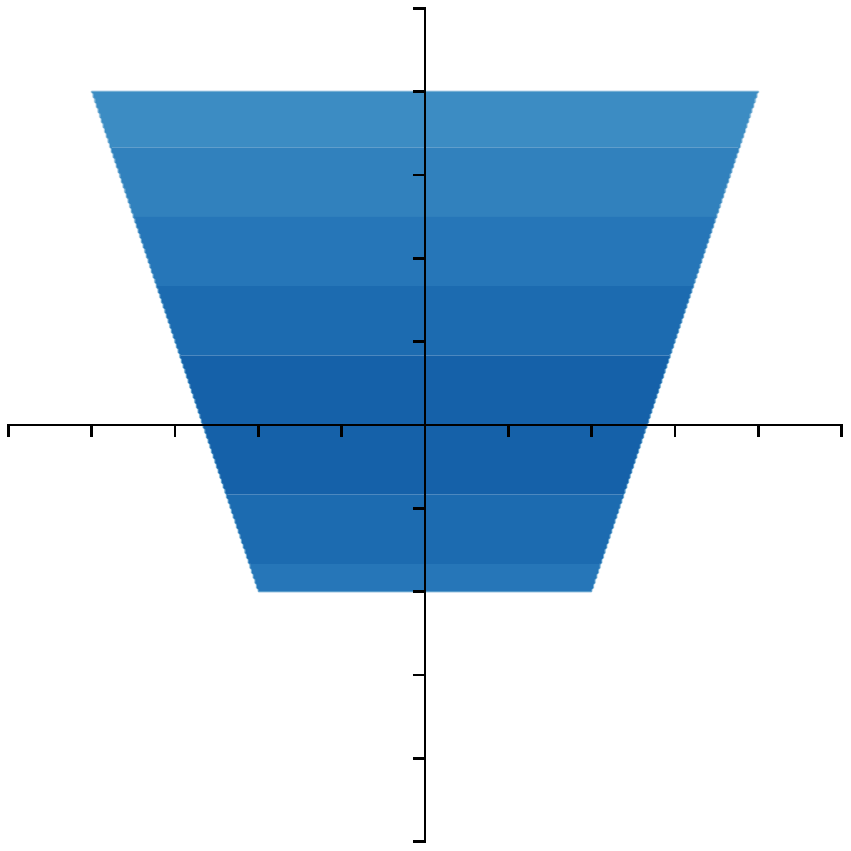}
							& \includegraphics[width=1.6cm]{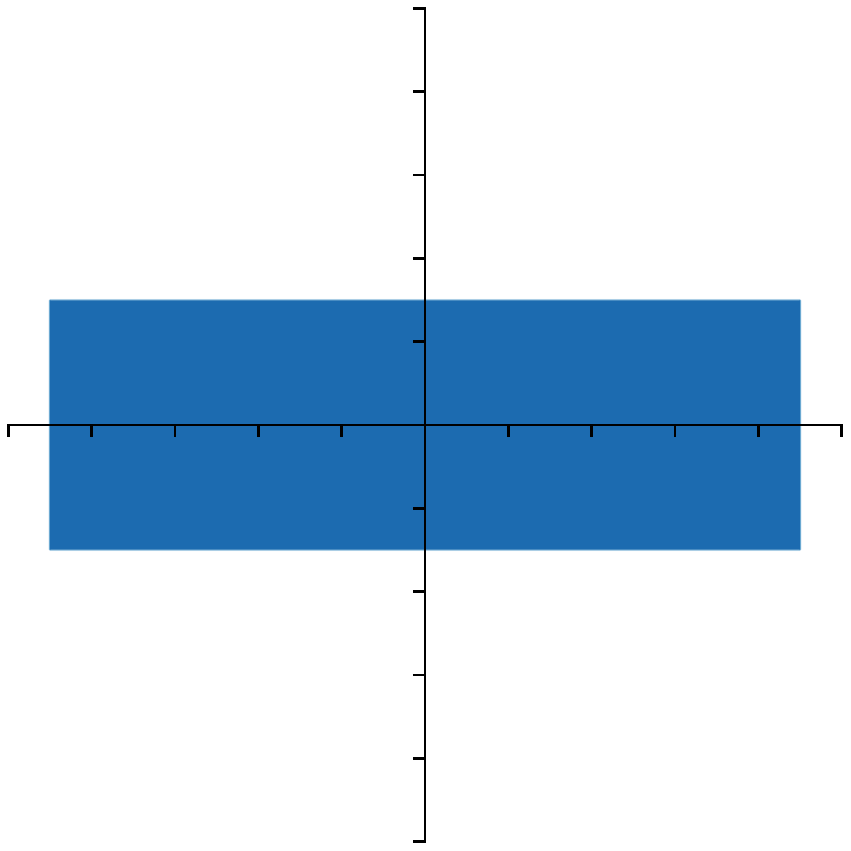}
							& \includegraphics[width=1.6cm]{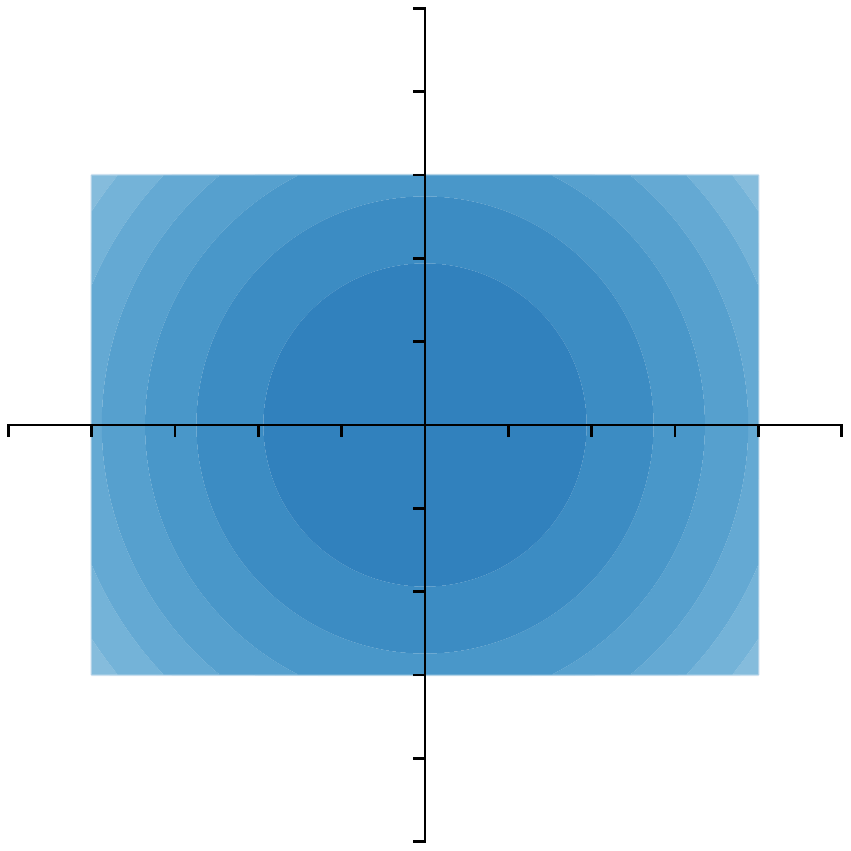}
							& \includegraphics[width=1.6cm]{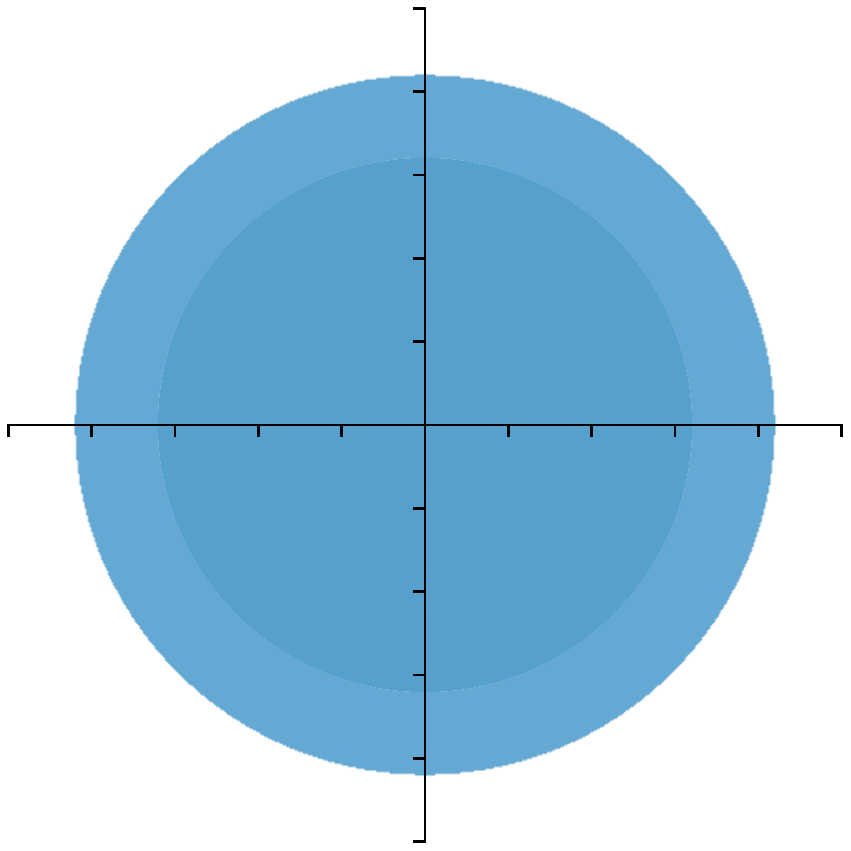} \\
		Sensor intensity $I$ & 316.91 & 937.76 & 800.24 & 641.25 & 1096.06 & 1428.25 \\
		\hline
	\end{tabular}
\end{table*}

For trajectory comparison, agent initial positions and directions are appointed according to given values, same for each control method, while for the search performance evaluation they are generated randomly in each run of Monte-Carlo simulation.

All considered methods are tested using Dubins motion model in all three test cases. In addition to Dubins model, the kinematic motion model is used in the first case in order to show possible variation of search performance due to the movement constraint in Dubins model. 
Although agent trajectories obtained with kinematic and Dubins motion models differ, search performance is not greatly influenced by trajectory curvature constraint in Dubins model for all control methods except RHC.
Keeping this in mind and the fact that real UAV motion is constrained by turning radius, the other two cases are tested only with Dubins motion model.

Plots of search agents trajectories for certain cases are presented in order to demonstrate diversity of considered control algorithms. Additionally, for the HEDAC controlled scenarios, the animated visualization of the search is provided in the supplementary material.

In order to determine and compare the search performance, 20 simulations are conducted for each scenario. This enables a Monte-Carlo method where sampling includes the stochastic effects of target and agent initial locations and the stochasticity of the target detection process.
Considering the average $E$ across 20 search simulation runs, the expected search time for ensuring 90\% of detection is determined. Since 90\% detection rate corresponds to $E=0.1$, the time $t_{90\%}$ needed to reach it is calculated, using simple root-finding method, in order to satisfy $E(t_{90\%}) = 0.1$. In author's opinion, the expected time needed to reach certain probability of finding target reveals intuitive measure of search success. The time is a critical resource in the real-world search campaigns.

\subsection{Test 1: Gaussian target distribution}

When estimating the search area, a very common case is to explore the region around a certain point, i.e. the last known location of the target.
The Gaussian function is a suitable choice for target occurrence probability density.
Similar estimation was used in the 2009 Air France airplane wreck search \cite{stone2014search}.
In this case we consider the 1000 m x 1000 m domain with 
\begin{equation}
	m_0(x,y) = A \cdot \exp \left(-\left({\frac {(x-x_{o})^{2}}{2\sigma _{x}^{2}}}+{\frac {(y-y_{o})^{2}}{2\sigma _{y}^{2}}}\right)\right)
	\label{eq:point_gaussian_distribution}
\end{equation}
where the center of the Gaussian function is $x_0=500$ m, $y_0=500$ m and the standard deviation $\sigma_x=\sigma_y=150$ m. The amplitude $A$ is calculated in order to satisfy \eqref{eq:occurence_normalization}.

The domain is discretized using $250 \times 250$ uniform rectangular grid.
In this test, the search is conducted using 5 identical mobile agents with parameters given in Table~\ref{tab:agent_parameters}.
The search simulation is performed for 600 s with time integration step $\Delta t = 0.25 \text{ s}$.
The HEDAC method parameters used for this case are: $\alpha = 0.03$ and $\beta = 4$.

\begin{figure*}[!ht]
	\centering
	\includegraphics[width=0.85\linewidth]{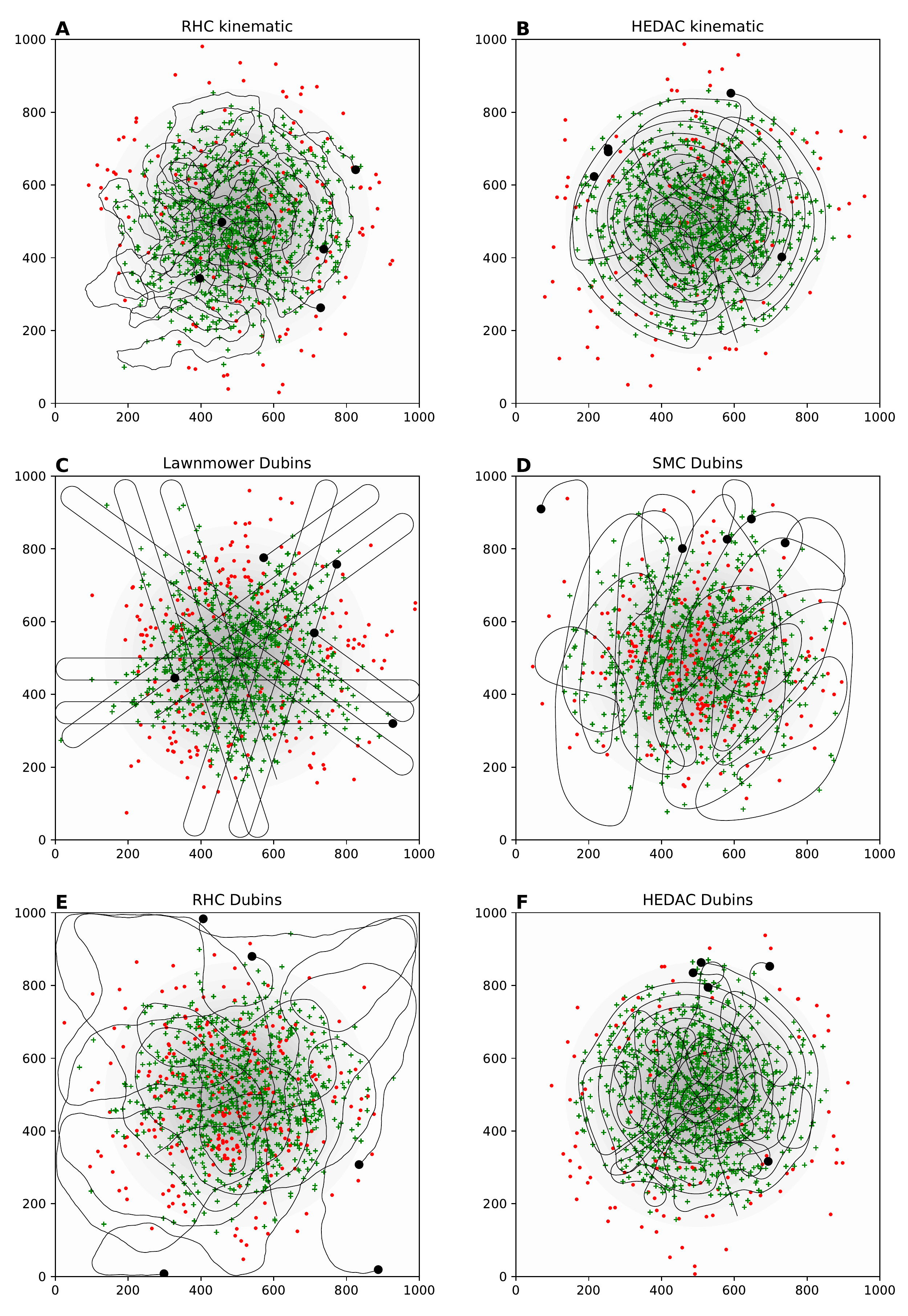}
	\caption{Test 1: Search agent trajectories comparison after 180 seconds of the search using kinematic motion model and RHC (A) and HEDAC (B) control method. The trajectories achieved using Dubins motion model after 180 seconds for Lawnmower (C), SMC (D), RHC (E) and HEDAC (F) search control. Dots and crosses (in online version of the manuscript colored red and green, respectively) mark undetected and detected targets, respectively. Agent locations at $t=180~\textmd{s}$ are denoted as black circles. (color figure available online)}
	\label{fig:test01_trajectories}
\end{figure*}

In Figure~\ref{fig:test01_trajectories}, agent trajectories, accomplished with all four considered methods for first 180 seconds of the search, are shown. The initial positions and directions of agents are defined as $x_{i,0}=500 + 70 \cdot i \cdot \cos((i-1)2\pi/5)~\text{[m]}$, $y_{i,0}=500 + 70 \cdot i \cdot \sin((i-1)2\pi/5~\text{[m]}$ and $\theta_{i,0}=(i-1)*\pi/5 + \pi$, for $i=1, ..., 5$. Only results for Dubins motion model are shown for Lawnmower and SMC algorithm, since there are no significant differences in search performance compared with kinematic model. For RHC and HEDAC methods, trajectories for both kinematic and Dubins motion models are presented in order to show difference in the movement of agents due to the turning radius constraints conditioned by the Dubins model. Observed agent control methods achieve fundamentally different motions. In contrast to geometrically regular Lawnmower trajectories, other control methods (SMC, RHC and HEDAC) produce chaotic behavior of search agents. RHC accomplished trajectories are stochastic due to underlaying PSO optimization method but SMC and HEDAC trajectories also seem rather stochastic, even though both methods are purely deterministic. If comparing trajectories for RHC, one can conclude that the Dubins constraint really aggravate
the optimization problem and makes it harder to achieve suitable paths. The minimal turning radius is perceptible in curvature of trajectories achieved using Dubins model.

The dynamics of search agents motion, controlled by HEDAC algorithm, can be observed in the Video 1 and Video 2 in supplementary material, for kinematic and Dubins motion model, respectively. Although the search is somewhat greedy, as all agents tend to go to regions of highest target probability, the cooperation of their motion is successful. 

\begin{figure*}[!ht]
	\centering
	\includegraphics[width=0.85\linewidth]{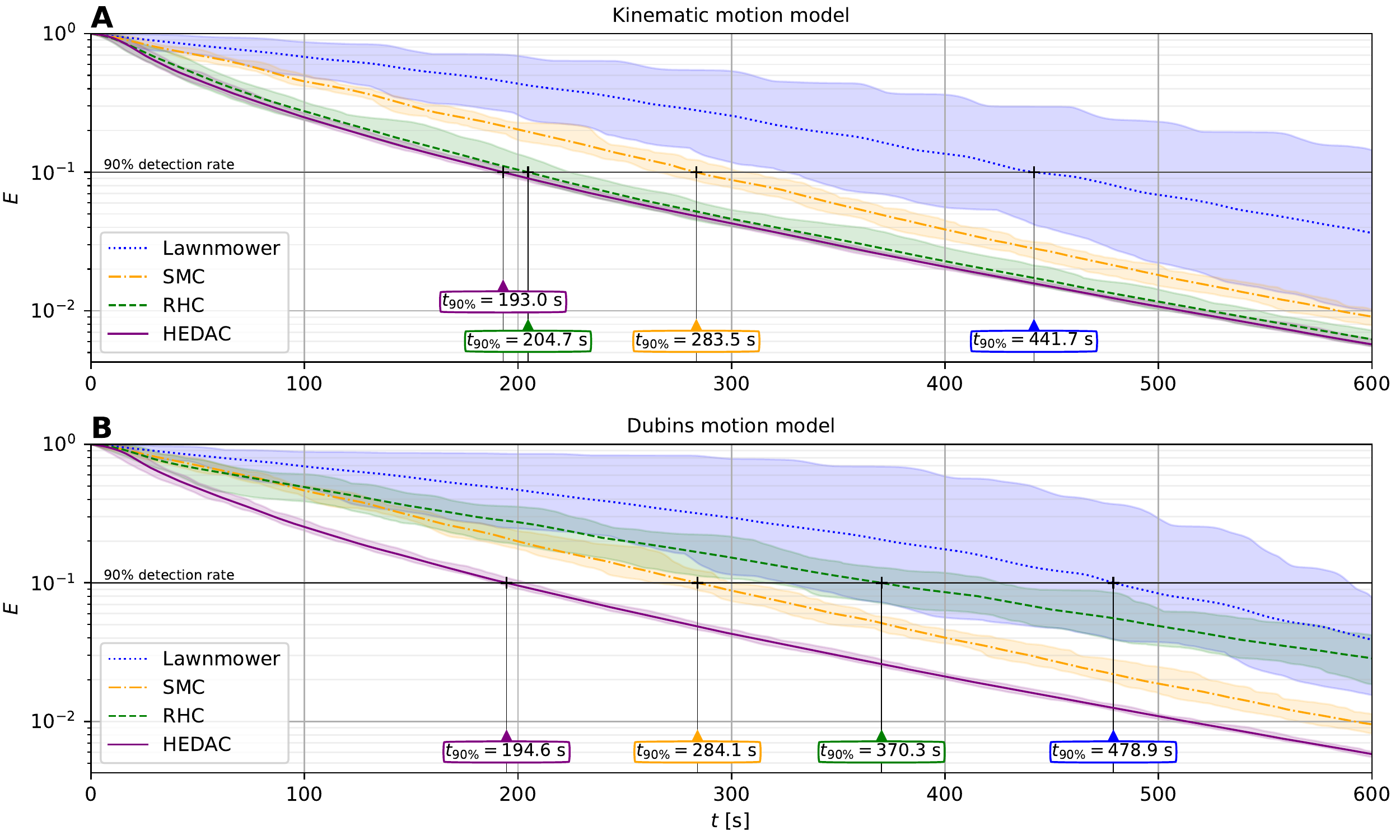}
	\caption{Test 1: Convergence comparison of Lawnmower, SMC, RHC and HEDAC search control methods for kinematic (A) and Dubins (B) motion model. The plot shows average value of $E$ (thick lines) and minimum and maximum values of $E$ (bounds of shaded areas) for each method over time, for 20 search simulation runs. Average time to reach 90\% detection rate, labeled as $t_{90\%}$, is shown for each scenario. (colored figure available online)} 
	\label{fig:test01_convergence}
\end{figure*}

The plot in Figure~\ref{fig:test01_convergence} shows convergence of target presence probability $E$ over time of the search process. The average values of $E$ across 20 conducted simulations indicate superiority of HEDAC over SMC and RHC, and especially Lawnmower algorithm. 
A small variation of $E$ indicates the robustness of HEDAC guided search in contrast to the alternative algorithms. HEDAC is able to achieve almost identical convergence of $E$ regardless of target positions, agent initial positions and detection stochasticity.

The search directed by HEDAC algorithm is clearly the best performing, regardless of the motion model used, with 193 s and 194.6 s on average to reach 90\% detection rate, for kinematic and Dubins motion model, respectively. 
The influence of Dubins motion constraint to search efficiency is clearly evident for RHC method. Although RHC with kinematic motion model ($t_{90\%}=204.7 s$) is slightly worse when compared with HEDAC, it needs almost double the time if Dubins constraint is used (370.3 s).
The SMC algorithm achieves the same detection rate in about 50\% more time than HEDAC (283.5 s and 284.1 s for kinematic and Dubins motion model, respectively), which is, in real-world situations, a significant deficiency. The search with Lawnmower is very inefficient, with over double time needed to reach 90\% detection rate (441.7 s and 478.9 s for kinematic and Dubins motion model, respectively), when compared to HEDAC controlled search.

\subsection{Test 2: Regionally uniform target distribution}

The multi-agent search on an island-like shape target probability distribution is demonstrated in Test 2. A simple construction, using union and subtraction of circular shapes, provides combined shape with uniform target occurrence probability density. Formally, the region of target distribution $\Omega_m \in \Omega$ can be defined as follows:
\begin{equation}
	\Omega_m = \bigcup_{k=1}^{n^+} \Omega^+_k - \bigcup_{k=1}^{n^-} \Omega^-_k
	\label{eq:test02_domain}
\end{equation}
where $\Omega^+$ is the region of minuend circle while $\Omega^-$ is the region of subtrahend circle. $n^+$ and $n^-$ denote number of minuend and subtrahend circles, respectively.

The initial target occurrence field is calculated as uniform over $\Omega_m$ region:
\begin{equation}
m_0(\mathbf{x}) = \begin{cases}
1 & \text{ if } \mathbf{x} \in \Omega_m \\
0 & \text{ otherwise}
\end{cases}
\end{equation}
and scaled to satisfy \eqref{eq:occurence_normalization}.

The combined shape is settled inside of a 3000 m $\times$ 3000 m domain which is divided into 600 segments in each direction ($\Delta x = \Delta y = 5~\text{m}$). The radius and center coordinates of constituent circles are given in supplementary material.
A heterogeneous set of 6 search agents is used: three pairs with different both motion and detection characteristics for each pair.
In HEDAC method scenario, agents are directed using parameters $\alpha = 0.03$ and $\beta = 2$. The time integration step used is $\Delta t = 0.5 \text{ s}$ and the total time of the simulated search is 1800 s.
Initial locations and directions of agents, for trajectories comparison, are (1000~m, 500~m, 0), (400~m, 1000~m, $\pi/6$), (1500~m, 1000~m, $\pi/3$), (1500~m, 2000~m, $\pi/2$), (2700~m, 2000~m, $2\pi/3$) and (2300~m, 2600~m, $5\pi/6$).

\begin{figure*}[!ht]
	\centering
	\includegraphics[width=0.85\linewidth]{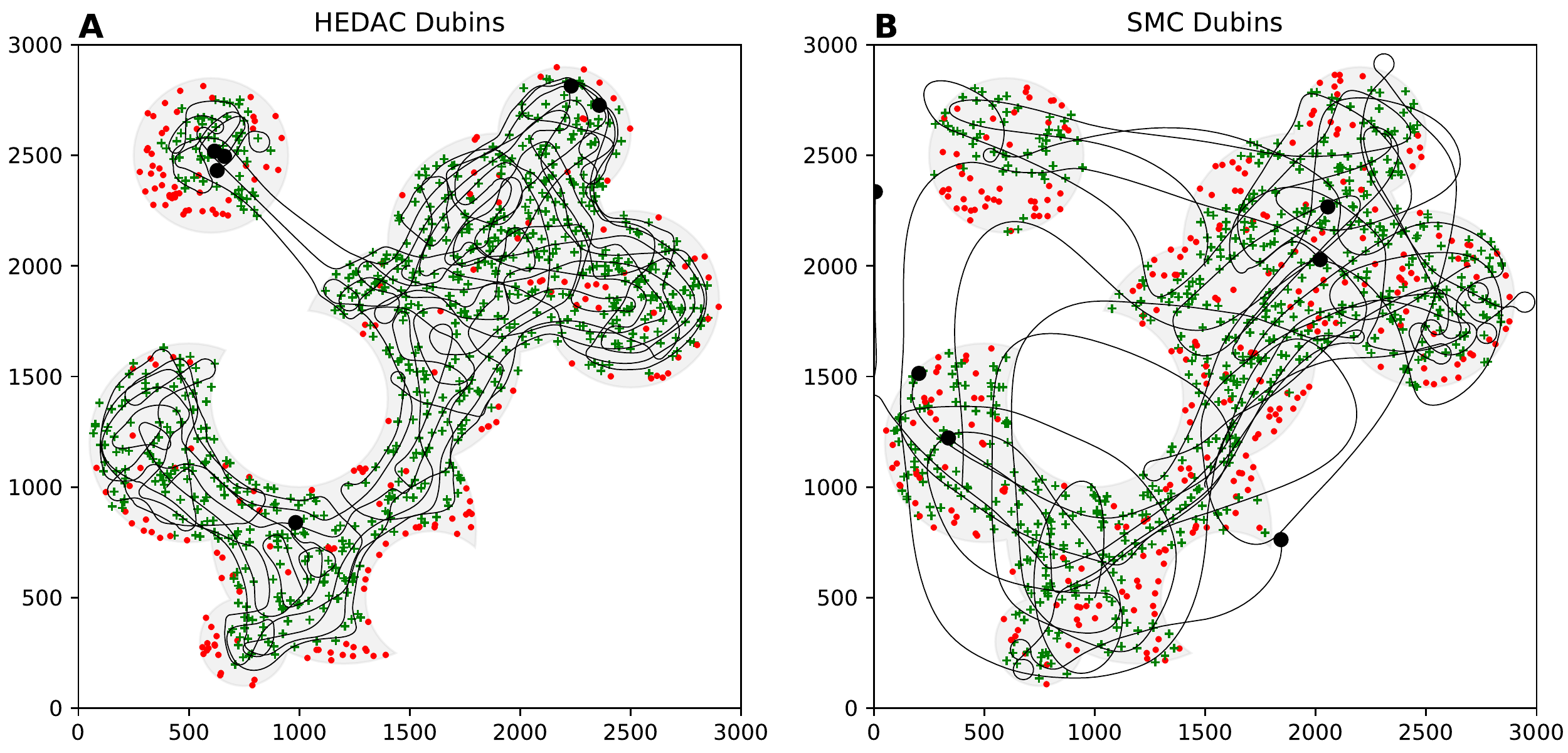}
	\caption{Test 2: Comparison of search agent trajectories after 600 s obtained using Dubins motion model with HEDAC (A) and SMC (B) search control. Dots and crosses mark undetected and detected targets, respectively. Current agent locations are marked as black circles. (color figure available online)}
	\label{fig:test02_trajectories}
\end{figure*}

The HEDAC drives agents within target occurrence shape with the exception of empty zone transit only to reach the isolated circle (Figure~\ref{fig:test02_trajectories}). Although the SMC also strives to cover the region of interest, the resulting trajectories are much broader and, consequentially, the search is not as efficient as with HEDAC. In SMC controlled scenario, there are many crossings between isolated regions, while HEDAC minimizes both the number of crossings and the length of paths between isolated regions.
Animation of the search performed using HEDAC algorithm can be found in Video 3 in supplementary material. One can easily notice a different agents motion and sensor characteristics, by observing agent velocity and turning radius, and the width and intensity of imprint in target probability field $m$ caused by detection sensor action.

\begin{figure*}[!ht]
	\centering
	\includegraphics[width=0.85\linewidth]{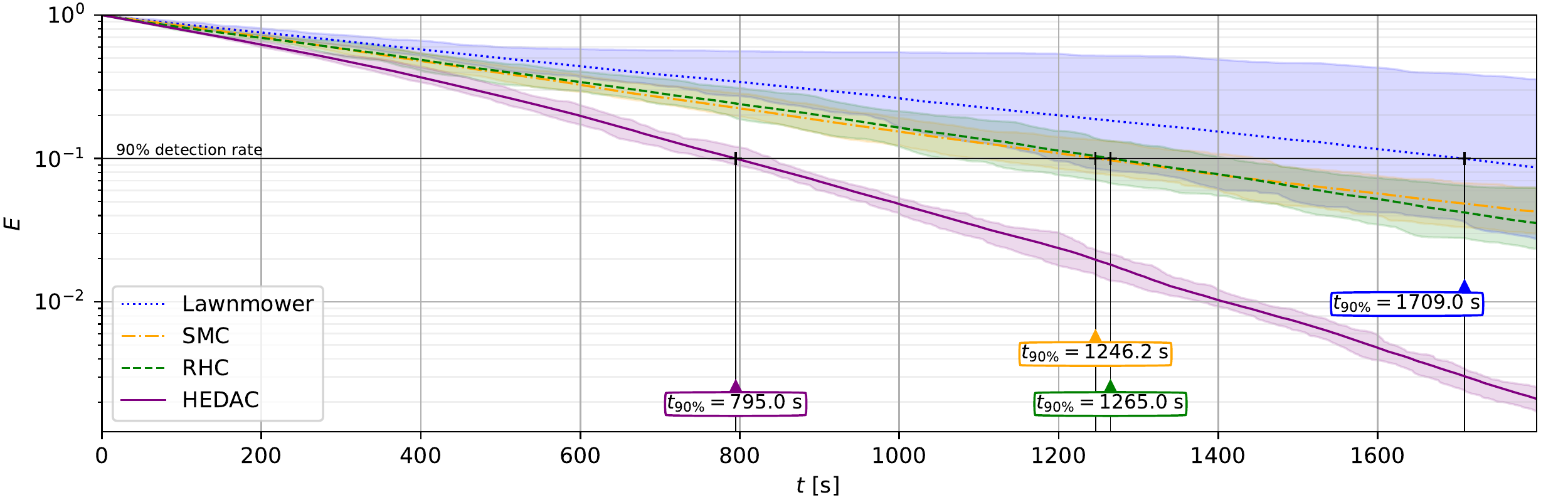}
	\caption{Test 2: Convergence comparison of Lawnmower, SMC, RHC and HEDAC search control methods using Dubins motion model. The plot shows average values of $E$ (thick lines) and minimum and maximum vlaues of $E$ (bounds of shaded areas) for 20 runs of each control algorithm. $t_{90\%}$ milestone is marked for each method. (color figure available online)}
	\label{fig:test02_convergence}
\end{figure*}

The convergence of $E$ for Test 2 is showed in Figure~\ref{fig:test02_convergence}, based on which one can reliably confirm the superiority of HEDAC over other compared control algorithms.

Due to more complicated shape of target distribution, the variation of total probability of undetected targets achieved with HEDAC is notably greater than in the first test case. Since the target distribution region is fragmented, it is inevitable for agents to pass over the "empty" region. Depending on the recurrence of such passing, the value of $E$ varies. The fluctuation of $E$ over 20 search simulations is stronger for the two alternative control algorithms, comapared with HEDAC.

Considering detection time (Figure~\ref{fig:test02_convergence}), once again the best results are obtained using HEDAC control algorithm where, on average, a 90\% detection is achieved in 795.0 s. Comparable to Test 1, the SMC underachieves HEDAC as measured with over 50\% longer search time (1246.2 s) to accomplish a 90\% detection rate. RHC struggles with overpassing unconnected regions and gives $t_{90\%}=1265~\text{s}$. Search inefficiency of the Lawnmower is evident by $t_{90\%}=1709~\text{s}$ which is roughly a double time needed than in HEDAC's search.

\subsection{Test 3: Nonuniform target distribution on a road network}

In the third test, a hypothetical, but potentially realistic, target search setting is considered. For a given arrangement of roads in the domain, the target occurrence probability density is determined so that presence of targets depends on the distance from the road. The probability follows the normal distribution law, similar to the formulation \eqref{eq:point_gaussian_distribution} used in Test 1. 
Using Gaussian radial basis function $\phi_\sigma$ over all road segments, the target occurrence density can be defined as follows:
\begin{equation}
m_0(\mathbf{x}) = A \sum_{s} \int_0^{L_s}  \phi_\sigma(\mathbf{x} - \mathbf{w}_s (l)) \diff l 
\label{eq:test03_target_distribution}
\end{equation}
where $\mathbf{w}_s (l)$ is position along road segment $s$ at distance $l$, $L$ is the length of the segment, and $A$ is constant adjusted so that $m_0$ satisfies \eqref{eq:occurence_normalization}. The broadness of target scattering is controlled with the standard deviation $\sigma$ of Gaussian radial basis function $\phi_\sigma$, which equals to 100 m.

A $4000~\text{m} \times 2000~\text{m}$ domain, discretized with $\Delta x = \Delta y = 5~\text{m}$,  is considered for this test case. The domain shape suites to demonstrate employment of the HEDAC method on rectangular domains, in contrast to square domain used in Test 1 and 2. The geometry of the road network, given with coordinates of straight road segments, is available in the supplementary material.
The search is performed for 3000 s using 5 agents whose characteristics are given in Table~\ref{tab:agent_parameters}. The time step is $\Delta t = 0.5 \text{ s}$, while the HEDAC parameters are: $\alpha=0.03$ and $\beta=2$.

\begin{figure}[!ht]
	\centering
	\includegraphics[width=\linewidth]{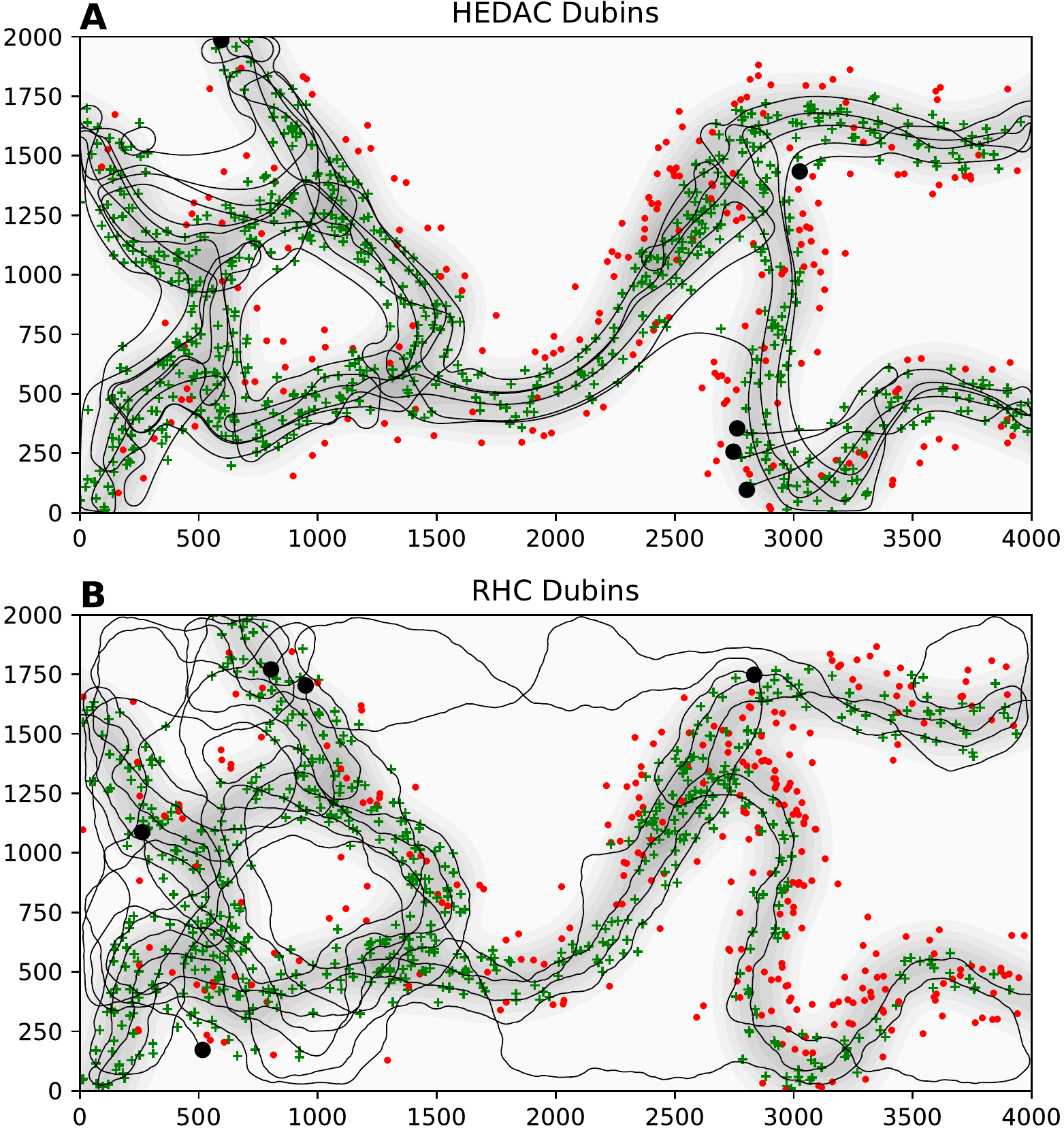}
	\caption{Search agent trajectories comparison for Test 3 after 600 s for HEDAC (A) and RHC (B) controlled search, both using Dubins motion model. Current positions of search agents are tagged with black circles. Undetected and detected targets are marked as dots and crosses, respectively. (color figure available online)}
	\label{fig:test03_trajectories}
\end{figure}

For agent trajectories comparison, agents are initially positioned at 5 points on the road at the boundary of the domain (where roads enter the domain). Agent initial directions correspond to inward normals of the domain boundary. 
Agents directed by HEDAC control accomplish trajectories which closely follow the region of road network and hence achieve reliable detection of scattered targets (Figure~\ref{fig:test03_trajectories}). Due to complicated shape of the target distribution, the global coverage tendency of SMC causes a more inefficient search. Although not under detailed investigation in this paper, the SMC showed more problems near the boundary of the domain. Due to limited turning in Dubins model, agents occasionally tend to go outside of domain boundaries but that is forcibly countered. Video 4 in the supplementary material shows the first 1800 s of HEDAC driven search.

\begin{figure*}[!ht]
	\centering
	\includegraphics[width=0.85\linewidth]{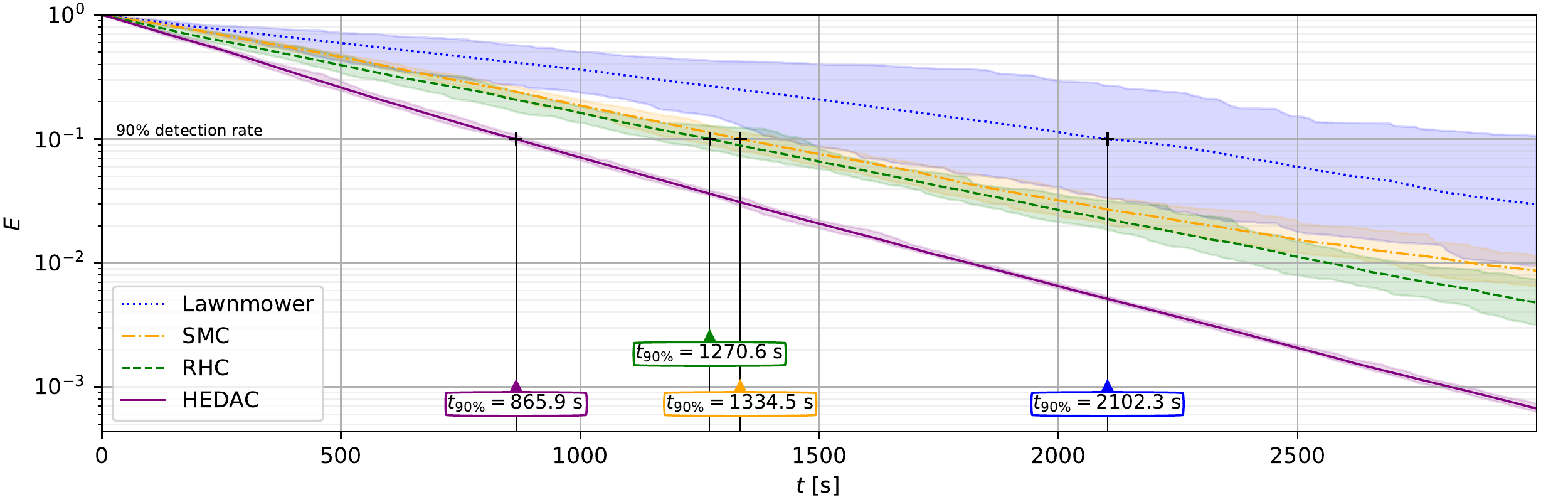}
	\caption{Test 3 convergence comparison of Lawnmower, SMC, RHC and HEDAC methods for Dubins motion model. The plot shows minimum, maximum and average value of $E$ for 20 runs of each control algorithm, and $t_{90\%}$ marks for each search method. (color figure available online)}
	\label{fig:test03_convergence}
\end{figure*}

Figure~\ref{fig:test03_convergence} shows the convergence of $E$ for all four compared control methods. Although the target distribution in this case is seemingly most complicated, due to interconnection of road network segments HEDAC performs a robust search with small deviation in convergence. The search success can be recognized in the time needed to accomplish 90\% detection probability. On average, HEDAC reaches this milestone in $t_{90\%}=865.9~\text{s}$, SMC in $t_{90\%}=1334.5~\text{s}$, RHC in $t_{90\%}=1270.6~\text{s}$ and Lawnmower in $t_{90\%}=2102.3~\text{s}$. 

\subsection{Computational efficiency}

All four considered search control algorithms differ in their approach and, consequentially, in their computational implementation and complexity.
In order to compare the methods from the aspect of computational efficiency, a simple benchmark is conducted using all three presented test cases.
The execution time of a single motion control step is measured for previously presented search scenarios. Various statistical computations, results saving and visualizations are excluded from the measurement. Finally, a mean execution time is calculated from first 100 search steps.

The benchmark is conducted on 64 bit PC, with 8 x Intel Core i7-4770 CPU @ 3.40 GHz and DDR3 32 GB 1600 MT/s with openSUSE Leap 15.1 operating system on Linux kernel 4.12.14. All four control methods are implemented using Python 3.7.3 and NumPy 1.16.4.

\begin{figure*}[!ht]
	\centering
	\includegraphics[width=0.85\linewidth]{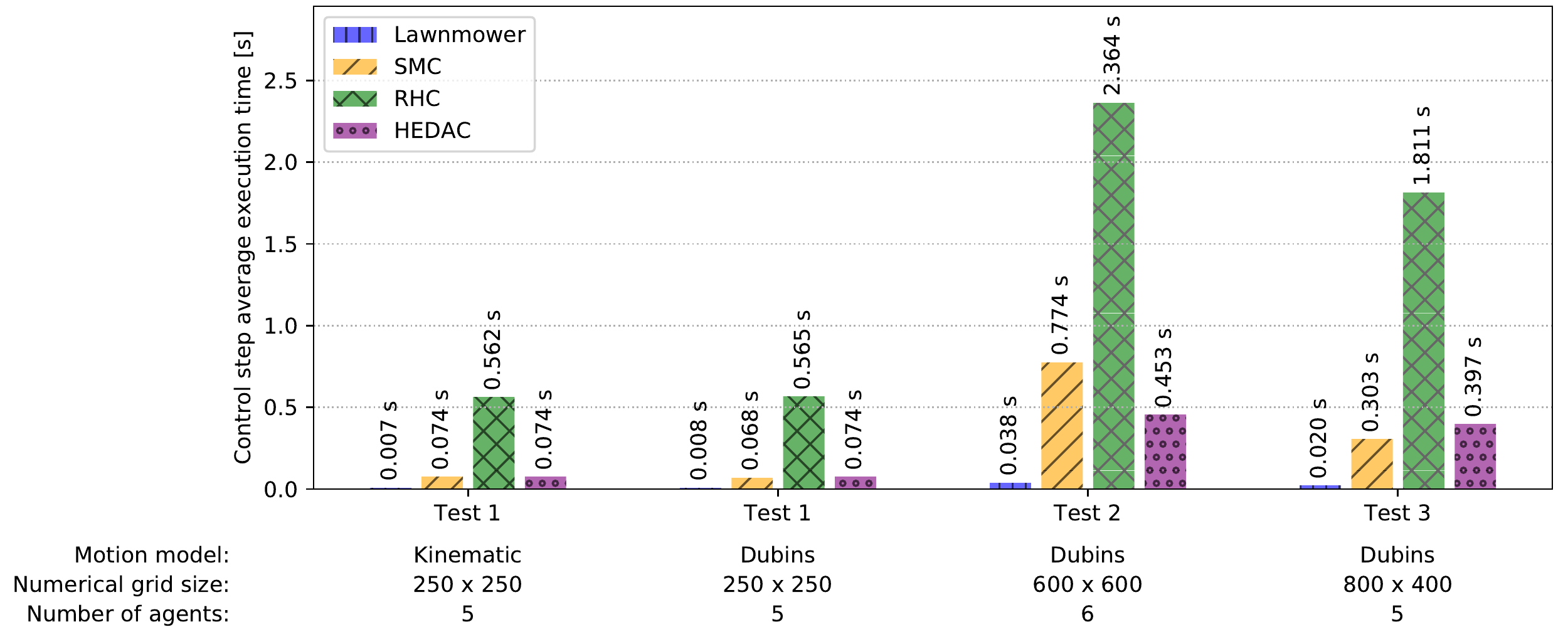}
	\caption{Benchmarking average control step execution time for all conducted scenarios. (color figure available online)}
	\label{fig:exec_benchmark}
\end{figure*}

The results of computational efficiency benchmark are summarized in Figure~\ref{fig:exec_benchmark}. The results shows no meaningful difference regarding used motion model for all four control methods, which is expected due to simplicity of Dubins constraint. Conventional Lawnmower method requires minimal calculation since entire path can be determined in advance of the search. Therefore, it is the absolutely fastest of the considered methods. The execution speed of SMC and HEDAC is comparable and with step time below one second, probably, employable for real world applications. An exception from the execution time vs. numerical grid size is visible for SMC in Test 2. In the author’s opinion this is caused by fast Fourier Transform (FFT) scalability where best performance are reached whit highest symmetry in the calculated terms which is when the size of input array is a power of 2. %2D FFT requires $\mathcal{O}(n_x n_y \log_2n_x \log_2n_y)$ operations
At last, the RHC is the slowest method in all observed scenarios. Due its underlaying optimization formulation, the RHC execution time not only depends on the numerical grid size but also on the number of search agents.

\subsection{HEDAC multi-agent search scalability}

One of the most important characteristics of multi-agent control is the scalability, i.e. taking the advantage of using multiple agents in order to speed-up the search. A search scenario defined in Test 1 is utilized for the HEDAC search efficiency analysis in which increasing number of search agents are used. All parameters, i.e. the domain and target distribution, agent motion and sensing characteristics and HEDAC parameters, are the same as in Test 1 and 20 search simulations are conducted for each $N=1,2,3,4,5,6,8,10,12,16,20$.

\begin{figure*}[!ht]
	\centering
	\includegraphics[width=0.85\linewidth]{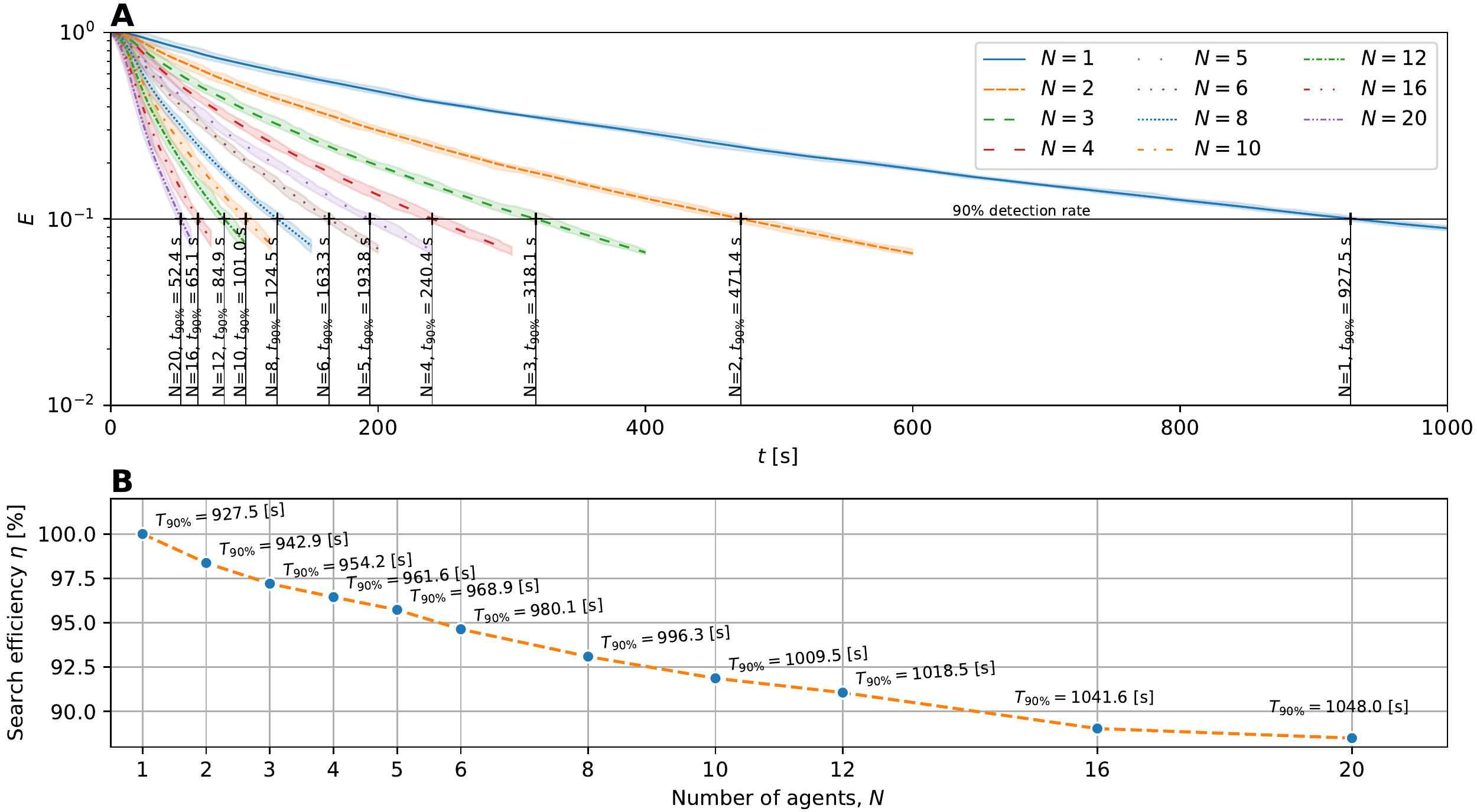}
	\caption{(A) The convergence of target presence probability $E$ for different number of agents $N$. (B) Analysis of the scalability by comparing the search efficiency for conducted searches with different number of agents used. (colored figure available online)}
	\label{fig:test04_scalability}
\end{figure*}

The convergence of average $E$ is, as expected, much faster when more agents are used as indicated by steeper curves in Figure~\ref{fig:test04_scalability}A as number of agents $N$ increases. Similarly to previous test cases, the search efficiency can be evaluated by $t_{90\%}$ - a time needed to achieve 90\% target detection rate. However, in order to provide a trustworthy and measurable evaluation of the search scalability, one needs to consider cumulative time spent by all agents which is needed to achieve the targeted 90\% detection rate. This scalability grade can be defined as $T_{90\%}(N) = t_{90\%}(N) \cdot N$. Finally, the search efficiency is defined relatively to single agent search as  $\eta(N) = T_{90\%}(1) / T_{90\%}(N)$.

The search efficiency gradually drops as number of used search agents increases (Figure~\ref{fig:test04_scalability}B). The effect of "over-searching", i.e. repeated passing over the area where sensing has already been conducted, is much more evident when larger number of search agents are used. Still, employing more agents with HEDAC control brings crucial advantage in terms of faster search without significantly deteriorating the search efficiency.

\section{Conclusion}

A success of the proposed methodology for multi-agent search relies on two important underlying presumptions. The first one is the precise estimation of target occurrence probability and the second is the ability of HEDAC multi-agent control algorithm to perform successful search according to estimated target distribution.   

In order to ensure trustworthy target occurrence distribution, a  probabilistic model is presented which incorporates the initial target distribution and search agent sensing action. This formulation includes the uncertainty of target detection based on search agent trajectories, and spatial and temporal characteristics of the sensor model. Described theory has a rigorous mathematical foundation and it can be potentially adjusted for various sensor models.

Heat Equation Driven Area Coverage (HEDAC) is adapted so that it guides agents according to the calculated probability (of undetected targets occurrence) density. A potential field is obtained with suitably designed stationary heat equation, where the source of the heat is distributed equivalently to the undetected targets probability field. Following the gradient of the potential, agents move to regions of undetected targets which maximizes the detection rate.
The proposed algorithm allows the use of heterogeneous search agents, considering both sensing/detection and motion characteristics. 

In order to test and validate the proposed methodology, the simulations of the search with Lawnmower, SMC, RHC and HEDAC control algorithms are performed on three hypothetical target search scenarios.
The results show significant advantage of HEDAC when compared with considered alternative methods. 
Use of simple search agent path generation, such as Lawnmower algorithm, can result in more than twice longer search for achieving the same success as with HEDAC.
SMC offers interesting spectral smoothing approach, but results indicate that HEDAC performs noticeably better for uncertainty target search. Although it can be tuned in order to achieve better search efficiency, the RHC method has severe drawbacks due to its optimization based approach. Satisfactory global search can be achieved with RHC only by setting up high-dimensional optimization problem which makes it unusable for real-world applications.
Expected time $t_{90\%}$ needed to realize 90\% detection rate is used as a measure to assess the search efficiency. Obtained results confirm the relation in search efficiency ($t_{90\%}$) between compared methods is roughly: HEDAC / SMC / Lawnmower = 1 / 1.5 / 2, while RHC achieves result in between other methods depending on given optimization/computational resources. In terms of computational efficiency, only RHC fails to attain real-time search control.

Analysis of HEDAC scalability showed that the proposed multi-agent motion control can be successfully used with large number of search agents. The decline in the search efficiency is easily compensated by the benefit of faster exploration when more agents are employed.

The proposed heterogeneous multi-agent control algorithm for target search in uncertain conditions offers a great flexibility and robustness. Beside initial target occurrence estimation, no further preprocessing and preparation is needed. The algorithm can handle any distribution of target occurrence and any number of search agents with different motion and sensor properties. 

Many improvements of the proposed search methodology are possible and here we briefly provide some ideas for future research on this topic.
Although the proposed method relies on centralized control, possibly it can be extended to fully autonomous decentralized multi-agent control. 
Further improvements are possible if more complex and comprehensive agent movement models are considered, for example a realistic multi-rotor UAV motion model. It should be noted that any agent motion model which can be controlled by a heading vector should be able to cooperate with proposed search system. Adding and removing agents from the search as it goes by is another possible enhancement of proposed framework. One of the most exciting enrichments of the method is the application in uncertainty search in unsteady conditions such as search for debris/people drifting at sea. The utilized target detection model is very simple and one should use a more accurate and detailed detection model, customized both temporally and spatially for used sensing equipment and image processing algorithms, in order to maximize the search success. A real experiment with adequately equipped UAV's for real-world search would be a further strong confirmation of the proposed methodology. Even though the search convergence of HEDAC multi-agent control method is shown in computational experiments, it still lacks rigorous mathematical proof of such convergence which requires a deeper theoretical study of the proposed feedback control model. 

Author believes the presented work serves as a significant advancement in a multi-agent search in uncertain conditions, which can potentially bring efficient real-world applications closer to reality.   

\section*{Acknowledgments}

This research is supported in part by the University of Rijeka under the project number 17.10.2.1.04

\bibliographystyle{plain}
\bibliography{bibliography}

\end{document}